\documentclass{article}

\usepackage{arxiv}
\usepackage{amssymb}
\usepackage{amsmath}
\usepackage{empheq}
\usepackage{amscd}
\usepackage{graphicx}
\usepackage{graphics}
\usepackage[noadjust]{cite}
\usepackage{amsthm}
\usepackage{caption}
\usepackage{multirow}
\usepackage{array}
\usepackage{rotating}

\usepackage{indentfirst}
\usepackage{tikz}
\usepackage{calc}
\usepackage{epsfig}
\usepackage{natbib}
\usepackage[makeroom]{cancel}
\usepackage{multicol}
\usepackage{csquotes}
\usepackage{algorithm}
\usepackage{algorithmic}
\usepackage{enumitem}
\usepackage{hyperref}
\usepackage{url}
\usepackage{bbm}  
\hypersetup{colorlinks,linkcolor={blue},citecolor={blue},urlcolor={blue}}

\usepackage{timet}
\usepackage{graphicx}
\usepackage{multirow}
\usepackage{multicol}
\usepackage{lipsum}
\usepackage{timet}
\usepackage{epsfig}
\usepackage{amsmath}
\usepackage{amsfonts}
\usepackage{amssymb}
\usepackage{color}
\numberwithin{equation}{section}
\usepackage{caption}
\usepackage{float}
\usepackage{subcaption}
\usepackage{graphics}
\usepackage{xcolor,graphicx}
\usepackage{csquotes}
\usepackage{mathtools}
\usepackage{optidef}

\usepackage{multicol}
\usepackage{empheq}
\usepackage{graphicx}
\usepackage{amssymb}
\usepackage{amsmath}
\usepackage{caption}
\usepackage{float}
\usepackage{subcaption}
\usepackage{graphicx}
\usepackage{hyperref}
\hypersetup{colorlinks,linkcolor={blue},citecolor={blue},urlcolor={blue}}
\usepackage{bbm}  
\newcommand{\vertiii}[1]{{\left\vert\kern-0.25ex\left\vert\kern-0.25ex\left\vert #1 
    \right\vert\kern-0.25ex\right\vert\kern-0.25ex\right\vert}}
\DeclareMathOperator*{\minimize}{minimize}

\newcommand{\eref}[1]{equation \ref{#1}}                  
\newcommand{\fref}[1]{Figure \ref{#1}} 
\begin{document}
\title{On the connection between  WRI and FWI: Analysis of the nonlinear term in the Hessian matrix}

\author{  \href{https://orcid.org/0000-0002-9879-2944}{\includegraphics[scale=0.06]{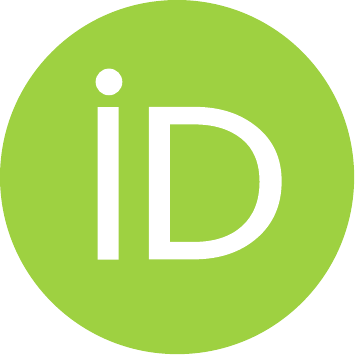}\hspace{1mm}Ali Gholami} \\
  Institute of Geophysics, University of Tehran, Tehran, Iran.
  \texttt{agholami@ut.ac.ir} \\ 
  Institute of Geophysics, Polish Academy of Sciences, Warsaw, Poland.  \texttt{agholami@igf.edu.pl} \\  
  \And
  \href{http://orcid.org/0000-0003-1805-1132}{\includegraphics[scale=0.06]{orcid.pdf}\hspace{1mm}Hossein S. Aghamiry} \\
  University Cote d'Azur - CNRS - IRD - OCA, Geoazur, Valbonne, France. 
  \texttt{aghamiry@geoazur.unice.fr}
\And
\href{http://orcid.org/0000-0002-4981-4967}{\includegraphics[scale=0.06]{orcid.pdf}\hspace{1mm}St\'ephane Operto} \\ 
  University Cote d'Azur - CNRS - IRD - OCA, Geoazur, Valbonne, France. 
  \texttt{operto@geoazur.unice.fr}
  }

\renewcommand{\shorttitle}{On the connection between  WRI and FWI, Gholami et al.}

\maketitle

\begin{abstract}
Implementation of the standard full waveform inversion (FWI) poses difficulties as the initial model offsets from the true model.
The wavefield reconstruction inversion (WRI) was proposed to mitigate these difficulties by relaxing the wave-equation constraint. In this abstract, working on the nonlinear term in the Hessian matrix of FWI, we develop a new approximate Hessian as an Augmented Gauss-Newton (AGN) Hessian including second-order derivative information. Moreover, we establish an intimate connection between an updating formula which results from approximate solve of the Newton's method with the AGN Hessian on the FWI problem and the WRI method. Our analysis opens new perspectives for developing efficient algorithms for FWI based on the Newton's method and highlights the importance of the nonlinear term in the Hessian matrix, which is ignored in most cases. 
\end{abstract}

\section{Introduction}
Full waveform inversion (FWI) has been widely accepted as an accurate method for the computation and characterization of subsurface model parameters by matching predicted to observed seismograms \citep{Tarantola_1988_TBI,Pratt_1998_GNF,Virieux_2009_OFW}. 
The Newton's method is an ideal option for solving the FWI problem but, in practice,  due the big size of the model parameters and the data to be inverted, the computational cost of implementing the Newton's method with full potential can be prohibitively high. Furthermore, in the case of poor initial models, the Hessian matrix can be ill-conditioned and indefinite, which causes an additional difficulty in determining an appropriate Newton step as an approximate solution of the Newton system of equations via the conjugate gradient (CG) method. 
These difficulties make the most practical algorithms to relay on simplified approximations of the Hessian matrix, such as the preconditioned steepest descent method.
However, when the initial model is of poor quality the associated data residual will be large, making the nonlinear term in the Hessian matrix, which is ignored in most cases,  important for convergence \citep{Pratt_1998_GNF}. 
Recently, attempts have been made to account for second-order information in the Hessian matrix such as truncated Newton method \citep{Metivier_2017_TRU}, l-BFGS and Anderson acceleration \citep{Yang_2021_AAS,Aghazade_2022_AAA}; however, the role of the nonlinear part of the Hessian is still not well understood.

Alternative methods have been proposed to increase the robustness of FWI with respect to the initial model. 
The wavefield reconstruction inversion (WRI) \citep{VanLeeuwen_2013_MLM}  increases the robustness by relaxing the wave-equation constraint. 
This wave-equation relaxation, or source extension, brings advantages for the inversion such as stability and robustness with respect to the initial model.
Recently, the analogy between the FWI misfit function and the WRI function was highlighted through a covariance matrix in the data space \citep{vanLeeuwen_2019_ANO,Symes_2020_WRI,Gholami_2022_EFW}.
Introducing this covariance matrix in the data space may be a leverage to move some effects of the Hessian in the gradient \citep{Thierry_1999_FRB,Metivier_2015_ASY}.

In this abstract, we reveal another connection between WRI and FWI. 
We explore the structure of the nonlinear term in the FWI Hessian matrix involving second-order derivative information and build a new approximation of this term as a weighted Gauss-Newton (GN) Hessian that leads to an augmented GN (AGN) approximation of the Hessian. This structured Hessian is then implemented efficiently and approximately to update the model by a sequential approach based on factorization. We clarify the connection between this second-order update formula, which results from an approximate solve of the Newton's method with the AGN Hessian on the FWI problem, and the update formula provided by the WRI. 
%
%
\section{Theory}
First, we propose a new approximate Hessian matrix for FWI by augmenting the standard GN Hessian with a term that involves second-order derivative information and provide an approximate while efficient update formula based on the new approximate Hessian. Then, we explore the mathematics of the WRI and determine its intimate connection with our updating formula, which results from an approximate Newton's method on the FWI problem.
\section{FWI}
The FWI problem is expressed as the nonlinear least-squares (LS) problem \citep{Tarantola_1988_TBI,Pratt_1998_GNF}.
\begin{equation} \label{FWI}
\minimize_{\bold{m}}~\frac{1}{2}\|\bold{d}^{*}-\bold{P}\bold{A(m)}^{-1}\bold{b}^{\!*}\|_2^2,
\end{equation}
where $\bold{m}$ is the model parameters, 
$\bold{P}$ is the sampling operator,
$\bold{A(m)}$ is the PDE operator, 
$\bold{b}^{\!*}$ is the source, 
and $\bold{d}^{*}$ is the recorded seismic data.
In this paper, we limit ourselves to the scalar acoustic wave equation with constant density, namely, $\bold{A(m)} = \bold{m} \circ \partial_{tt} - \boldsymbol{\Delta}$, where $\circ$ is the element-wise product operator.
For sake of compactness, we consider one source.

Standard methods for solving the nonlinear objective function \eqref{FWI} are gradient based in which the parameters are updated at iteration $k$ as $\bold{m}_{k+1} = \bold{m}_k + \delta\bold{m}_k$,
where the model update $\delta\bold{m}_k$ is calculated by solving the following Newton system:
\begin{equation} \label{NewtonSystem}
\bold{H}_k \delta\bold{m}_k=-\bold{g}_k,
\end{equation}
and $\bold{H}_k\equiv \bold{H(m}_k)$ and $\bold{g}_k\equiv \bold{g(m}_k)$ are the Hessian matrix and the gradient vector, which are explored in what follows.

Let us define $\delta\bold{d}(\bold{m})=\bold{d}^*-\bold{d(m)}$ and $\bold{d(m)}=\bold{PA(m})^{-1}\bold{b}^{\!*}$,
then the gradient $\bold{g}$ of the misfit function to be minimized, eq.~\ref{FWI}, is given by
\begin{equation} \label{A1}
\bold{g}=\nabla_{\!\bold{m}}\left(\frac{1}{2}\|\delta\bold{d}\|_2^2 \right)= \left(\nabla_{\!\bold{m}} \delta\bold{d}\right)^T\delta\bold{d}=\bold{J}^T\delta\bold{d},
\end{equation}
where $\nabla_{\!\bold{m}}$ is the gradient with respect to $\bold{m}$, superscript $T$ is the matrix transposition, and $\bold{J}$ is the first partial derivative (Jacobian) matrix of the data residual defined as 
\begin{equation} \label{J}
\bold{J}=\bold{S}\bold{L},
\end{equation}
with
\begin{equation} \label{Virsource}
\bold{S}=\bold{PA}(\bold{m})^{-1},~~~ \bold{L}=\frac{\partial \bold{A(m)}}{\partial \bold{m}}\bold{u}, ~~~\bold{u}=\bold{A}(\bold{m})^{-1}\bold{b}^{\!*}.
\end{equation}
In the above equations and what follows, the variables $\bold{d}$, $\delta\bold{d}$, $\bold{A}$, $\bold{S}$, $\bold{L}$, $\bold{J}$ and $\bold{u}$ depend on $\bold{m}$.

The $ij$th component of the Hessian matrix is
\begin{align} \label{A2}
\bold{H}_{ij}&= \underbrace{\left(\frac{\partial \bold{d}}{\partial \bold{m}_i}\right)^T  \left(\frac{\partial \bold{d}}{\partial \bold{m}_j}\right)}_{[\bold{J}^T\bold{J}]_{ij}}
-  \underbrace{\left(\frac{\partial^2 \bold{d}}{\partial\bold{m}_j\partial \bold{m}_i}\right)^T\delta\bold{d}}_{\bold{R}_{ij}},
\end{align}
where $\bold{J}^T\bold{J}$ is the GN term of the Hessian, which includes the first-order partial derivatives, and $\bold{R}$ includes the second-order partial derivatives \citep{Pratt_1998_GNF}.
In order to explore the structure of the nonlinear term, $\bold{R}$, we begin by the observation equation $\bold{d}=\bold{Pu}$, from which we get
\begin{align} \label{diff2_d}
\frac{\partial^2 \bold{d}}{\partial\bold{m}_j\partial \bold{m}_i}=
\bold{P}\frac{\partial^2 \bold{u}}{\partial\bold{m}_j\partial \bold{m}_i}.
\end{align}
We differentiate both sides of the wave equation, $\bold{A}\bold{u}=\bold{b}^*$, with respect to $\bold{m}_i$, getting
\begin{equation} \label{diff_wq}
\frac{\partial \bold{A}}{\partial \bold{m}_i}\bold{u}+ \bold{A} \frac{\partial \bold{u}}{\partial \bold{m}_i}=\bold{0},
\end{equation}
and then differentiate both sides of \eref{diff_wq} with respect to $\bold{m}_j$, giving 
\begin{equation} \label{diff2_wq}
\frac{\partial^2 \bold{A}}{\partial \bold{m}_j \partial \bold{m}_i}\bold{u}+ \frac{\partial \bold{A}}{\partial \bold{m}_i} \frac{\partial \bold{u}}{\partial \bold{m}_j} + \frac{\partial \bold{A}}{\partial \bold{m}_j} \frac{\partial \bold{u}}{\partial \bold{m}_i}+\bold{A} \frac{\partial^2 \bold{u}}{\partial \bold{m}_j \partial \bold{m}_i}=\bold{0}, \nonumber
\end{equation}
from which, we obtain
\begin{equation} \label{diff2_u}
\frac{\partial^2 \bold{u}}{\partial \bold{m}_j\partial \bold{m}_i} = -\bold{A}^{-1}\left(
  \frac{\partial \bold{A}}{\partial \bold{m}_i} \frac{\partial \bold{u}}{\partial \bold{m}_j}+
   \frac{\partial \bold{A}}{\partial \bold{m}_j} \frac{\partial \bold{u}}{\partial \bold{m}_i}+
  \frac{\partial^2 \bold{A}}{\partial \bold{m}_j \partial \bold{m}_i} \bold{u}\right). \nonumber
\end{equation}
We plug $\frac{\partial \bold{u}}{\partial \bold{m}_k}$ from \eref{diff_wq} into the above equation and then the result into \eref{diff2_d}, obtaining
\begin{multline}\label{R}
\hspace{-0.25cm}
\bold{R}_{ij}
  =  -\left(\frac{\partial \bold{A}}{\partial \bold{m}_j}\bold{u}\right)^T \bold{A}^{-T}\frac{\partial \bold{A}^T}{\partial \bold{m}_i}\bold{S}^T\delta\bold{d} 
  - \left(\frac{\partial \bold{A}}{\partial \bold{m}_i}\bold{u}\right)^T \bold{A}^{-T}\frac{\partial \bold{A}^T}{\partial \bold{m}_j}\bold{S}^T\delta\bold{d}  
 + \left(\frac{\partial^2 \bold{A}}{\partial \bold{m}_j \partial \bold{m}_i} \bold{u} \right)^T\bold{S}^T\delta\bold{d}. 
\end{multline}
In order to further simplify this expression, we explore the structure of $\bold{A}^{-T}\frac{\partial \bold{A}^T}{\partial \bold{m}_k}\bold{S}^T$.
We write
\begin{subequations}
\begin{align}
\frac{\partial (\bold{S}^T\bold{S})}{\partial \bold{m}_k}&=\frac{\partial \bold{S}^T}{\partial \bold{m}_k}\bold{S}+\bold{S}^T\frac{\partial \bold{S}}{\partial \bold{m}_k}\\
&= \frac{\partial \bold{A}^{-T}}{\partial \bold{m}_k}\bold{P}^T\bold{S} + \bold{S}^T\bold{P}\frac{\partial \bold{A}^{-1}}{\partial \bold{m}_k}\\
&= -\bold{A}^{-T}\frac{\partial \bold{A}}{\partial \bold{m}_k}\bold{A}^{-T}\bold{P}^T\bold{S} - \bold{S}^T\bold{P}\bold{A}^{-1}\frac{\partial \bold{A}}{\partial \bold{m}_k}\bold{A}^{-1} \label{key1}\\
&= -\bold{A}^{-T}\frac{\partial \bold{A}}{\partial \bold{m}_k}\bold{S}^T\bold{S} - \bold{S}^T\bold{S}\frac{\partial \bold{A}}{\partial \bold{m}_k}\bold{A}^{-1}, \label{key2}
\end{align}
\end{subequations}
where \eref{key1} is obtained by using the following identity for derivative of the inverse of a matrix: 
\begin{equation}
\frac{\partial \bold{A}^{-1}}{\partial \bold{m}}=-\bold{A}^{-1} \frac{\partial \bold{A}}{\partial \bold{m}} \bold{A}^{-1},
\end{equation}
then we right multiply both sides of \eref{key2} by $\bold{S}^T(\bold{S}\bold{S}^T)^{-1}$ and  rearrange the terms, getting
\begin{align}
\bold{A}^{-T}\frac{\partial \bold{A}^T}{\partial \bold{m}_k} \bold{S}^T
&=  -(\bold{S}^T\bold{S})\frac{\partial \bold{A}}{\partial \bold{m}_k}\bold{A}^{-1}\bold{S}^T(\bold{S}\bold{S}^T)^{-1} - \frac{\partial (\bold{S}^T\bold{S})}{\partial \bold{m}_k}\bold{S}^T(\bold{S}\bold{S}^T)^{-1}. \notag 
\end{align}
Substituting this equation into \eref{R} gives
\begin{multline}\label{RR}
\bold{R}_{ij}
   = \underbrace{\left(\frac{\partial \bold{A}}{\partial \bold{m}_i}\bold{u}\right)^T (\bold{S}^T\bold{S})\frac{\partial \bold{A}}{\partial \bold{m}_j}\delta\bold{u}}_{\bold{R}^{11}_{ij}}
+ \underbrace{\left(\frac{\partial \bold{A}}{\partial \bold{m}_i}\bold{u}\right)^T\frac{\partial (\bold{S}^T\bold{S})}{\partial \bold{m}_j}\delta \bold{b}}_{\bold{R}^{12}_{ij}} 
+ \underbrace{\left(\frac{\partial \bold{A}}{\partial \bold{m}_j}\bold{u}\right)^T (\bold{S}^T\bold{S})\frac{\partial \bold{A}}{\partial \bold{m}_i}\delta\bold{u}}_{\bold{R}^{21}_{ij}}
  + \underbrace{\left(\frac{\partial \bold{A}}{\partial \bold{m}_j}\bold{u}\right)^T\frac{\partial (\bold{S}^T\bold{S})}{\partial \bold{m}_i}\delta \bold{b}}_{\bold{R}^{22}_{ij}} \\
+ \left(\frac{\partial^2 \bold{A}}{\partial \bold{m}_j \partial \bold{m}_i} \bold{u} \right)^T\bold{S}^T\delta\bold{d}, \hspace{4cm}   
\end{multline}
where $\delta\bold{u}$ and $\delta\bold{b}$ are 
\begin{equation}
\label{delLub}
\delta\bold{u}=\bold{A(m)}^{-1}\delta\bold{b},~~~~~
\delta \bold{b} = \bold{S}^T(\bold{S}\bold{S}^T)^{-1}\delta\bold{d}. 
\end{equation}
%
\subsection{FWI with an Augmented Gauss-Newton Hessian}

Now we solve the Newton system in \eref{NewtonSystem} while using a new approximate Hessian obtained by using the Gauss-Newton Hessian and the $\bold{R}^{11}$ term of the nonlinear part.
We use $\bold{L}$ (equation \ref{Virsource}) and $\delta\bold{L} =\frac{\partial \bold{A}}{\partial \bold{m}}\delta\bold{u}$, then 
$\bold{R}\approx \bold{R}^{11}=\bold{L}^T \bold{S}^T\bold{S}\delta\bold{L}$ and device the following augmented GN (AGN) Hessian:
\begin{align} \label{EGN_app}
\bold{H} \approx  \bold{L}^T \bold{S}^T\bold{S}\bold{L} + \bold{L}^T \bold{S}^T\bold{S}\delta\bold{L},
\end{align}
from which the model update formula in \eqref{NewtonSystem} becomes
\begin{equation} \label{GN_update_factor}
\bold{L}_k^T\bold{S}_k^T\bold{S}_k(\bold{L}_k+\delta\bold{L}_k) \delta\bold{m}_k=-\bold{L}_k^T\bold{S}_k^T\delta{\bold{d}}_k.
\end{equation}
Now we explain how an approximate solve of this system leads to a model update that is exactly equal to that provided by the WRI method. 
Let us solve this equation by a sequential solve procedure similar to the procedure of solving a system of linear equations with (e.g. a lower-upper triangular) factorization.
One should note that sequential solve of the linear system \ref{GN_update_factor} leads to a different model update because the matrix $\bold{S}_k^T\bold{S}_k$ is rank deficient. \\
Consider that the virtual-source matrix $\bold{L}_k$ is full-rank, 
\begin{itemize}
\item[Step 1.] Multiply both sides of \eref{GN_update_factor} by $(\bold{L}_k^T)^{-1}$ to get
\begin{equation} \label{GN_update_factor2}
\bold{S}_k^T\bold{S}_k(\bold{L}_k+\delta\bold{L}_k) \delta\bold{m}_k=-\bold{S}_k^T\delta{\bold{d}}_k.
\end{equation}
Introduce the (scattering source) variable $\delta{\bold{b}}_k=-(\bold{L}_k+\delta\bold{L}_k)\delta\bold{m}_k$ and proceed with the following steps:
\item[Step 2.] Solve the equation $\bold{S}_k^T\bold{S}_k\delta{\bold{b}}_k=\bold{S}_k^T\delta{\bold{d}}_k$ for $\delta{\bold{b}}_k$, which can be considered as solving, in the LS sense, the linear system  $\bold{S}_k\delta{\bold{b}}_k=\delta{\bold{d}}_k$. After adding a damping term, we get
\begin{align} \label{deltab}
\delta{\bold{b}}_k = (\bold{S}_k^T\bold{S}_k+\epsilon \bold{I})^{-1}\bold{S}_k^T\delta{\bold{d}}_k
= \bold{S}_k^T(\bold{S}_k\bold{S}_k^T+\epsilon \bold{I})^{-1}\delta{\bold{d}}_k.
\end{align}
This damping is also applied to $\delta{\bold{b}}$ in \eref{delLub} which builds $\delta\bold{L}_k$.
This definition of $\delta\bold{b}_k$ makes the AGN Hessian depending on $\epsilon$. The limit case $\epsilon\to \infty$ results in $\delta\bold{b}_k=\bold{0}$, reducing the AGN Hessian to the standard GN Hessian. On the other hand, the limit case $\epsilon\to 0$ results in  $\delta\bold{b}_k$ with minimum energy that may be unstable due to the ill-conditioning of $\bold{S}_k\bold{S}_k^T$. 

\item[Step 3.] Solve the equation $(\bold{L}_k+\delta\bold{L}_k)\delta\bold{m}_k=-\delta{\bold{b}}_k$ for $\delta\bold{m}_k$, giving
\begin{equation} \label{GN_update_seq_data}
\delta\bold{m}_k=-\frac{\partial_{tt}(\bold{u}_k+\delta\bold{u}_k)\circ \delta{\bold{b}}_k}{|\partial_{tt}(\bold{u}_k+\delta\bold{u}_k)|^2}.
\end{equation}
\end{itemize}
An advantage of this sequential solve over directly solving \eref{GN_update_factor} for $\delta\bold{m}_k$ is that the data-domain Hessian matrix $\bold{S}_k\bold{S}_k^T + \epsilon \bold{I}$ in \eref{deltab} is source independent for fixed-spread acquisition and thus we need to construct it only once and use it for multiple sources.
However, due to the rank deficiency of matrix $\bold{S}_k^T\bold{S}_k$, this computational efficiency is achieved at the cost of approximately solving the original system. 

The numerator of \eref{GN_update_seq_data} is the zero-lag cross-correlation between the second partial derivative of the extended source wavefield, $\partial_{tt}(\bold{u}_k+\delta\bold{u}_k)$, and a receiver wavefield $\delta{\bold{b}}_k$, which is obtained by back propagation of a source term $(\bold{S}_k\bold{S}_k^T+\epsilon \bold{I})^{-1}\delta{\bold{d}}_k$, which is simply the data residual deconvolved by the data space Hessian, $\bold{S}_k\bold{S}_k^T+\epsilon \bold{I}$ \citep{Gholami_2022_EFW}. 

\section{WRI}
The WRI problem is obtained by relaxing the wave-equation constraint in \eref{FWI} and is expressed as \citep{VanLeeuwen_2013_MLM}
\begin{equation}  \label{WRI}
\minimize_{\bold{m},\bold{u}^e}~\|\bold{P}\bold{u}^e-\bold{d}^{*}\|_2^2+\mu \|\bold{A(m)}\bold{u}^e-\bold{b}^{\!*}\|_2^2,
\end{equation}
where $\bold{u}^e$ is the ``data-assimilated" wavefield, a wavefield that approximately satisfies the wave equation and the data simultaneously, $\mu$ is the penalty parameter that, for any given $\bold{m}$, balances between the data fit and the wave-equation fit.
This wave-equation relaxation brings advantages for the inversion such as stability and relative robustness with respect to the initial model.
The inversion begins with a small value assigned to $\mu$ such that simulated data closely match observed data with poor initial models. Then, the model is updated to reduce the wave-equation errors by properly tuning the value of $\mu$ \citep{Aghamiry_2019_IWR}. 
The objective function \eqref{WRI} is solved in an alternating mode.
Giving $\bold{m}_k$ at iteration $k$th, the associated wavefield, $\bold{u}_k^e$, is obtained by minimizing the objective function with respect to $\bold{u}^e$, giving
\begin{equation}
(\bold{P}^T\bold{P}+\mu \bold{A}_k^T\bold{A}_k)\bold{u}_k^e= \bold{P}^T\bold{d}^{*}+\mu \bold{A}_k^T\bold{b}^{\!*}.
\end{equation}
The solution $\bold{u}_k^e$ to this normal system also satisfies the following equation \citep{Gholami_2022_EFW}
\begin{equation} \label{ue}
\bold{A}_k\bold{u}_k^e= \bold{b}^{\!*} + \bold{S}_k^T(\bold{S}_k\bold{S}_k^T+\mu \bold{I})^{-1}\delta{\bold{d}}_k.
\end{equation}
For $\mu=\epsilon$, the second term at right hand side is exactly the scattering source in \eref{deltab}, which is a damped variant of the $\delta\bold{b}$ in  \eref{delLub} forming the scattered wavefield $\delta\bold{u}$. This clearly shows that $\bold{u}_k^e=\bold{u}_k+\delta\bold{u}_k$.
Having $\bold{u}_k^e$, WRI updates the model by  minimizing the objective function \eqref{WRI} with respect to $\bold{m}$, giving $\bold{m}_{k+1}$ as the minimizer of the source residual
\begin{align}  \label{WRIm}
\|\bold{A(m)}\bold{u}_k^e-\bold{b}^{\!*}\|_2^2=\|(\bold{m} \circ \partial_{tt}\bold{u}_k^e - \boldsymbol{\Delta}\bold{u}_k^e-\bold{b}^{\!*}\|_2^2,
\end{align}
that gives
\begin{equation} \label{WRIm2}
\bold{m}_{k+1} = \frac{\boldsymbol{\partial_{tt}\bold{u}_k^e\circ (\Delta}\bold{u}_k^e+\bold{b}^{\!*})}{|\partial_{tt}\bold{u}_k^e|^2}.
\end{equation} 
From \eref{ue} we have $\Delta\bold{u}_k^e+\bold{b}^{\!*}=\bold{m}_k\circ \partial_{tt}\bold{u}_k^e-\delta\bold{b}_k$ and thus \eqref{WRIm2} can be written as
\begin{equation} \label{WRIm3}
\bold{m}_{k+1} = \bold{m}_k - \frac{\partial_{tt}\bold{u}_k^e\circ \delta\bold{b}_k}{|\partial_{tt}\bold{u}_k^e|^2}.
\end{equation} 
Comparing this equation with \eref{GN_update_seq_data} and noting that $\bold{u}_k^e=\bold{u}_k+\delta\bold{u}_k$, we get that the model update provided by WRI can be obtained by FWI that approximately solves the Newton system with the AGN Hessian.
%
%
%
\begin{figure}[htb!]
\center
\includegraphics[width=1\columnwidth]{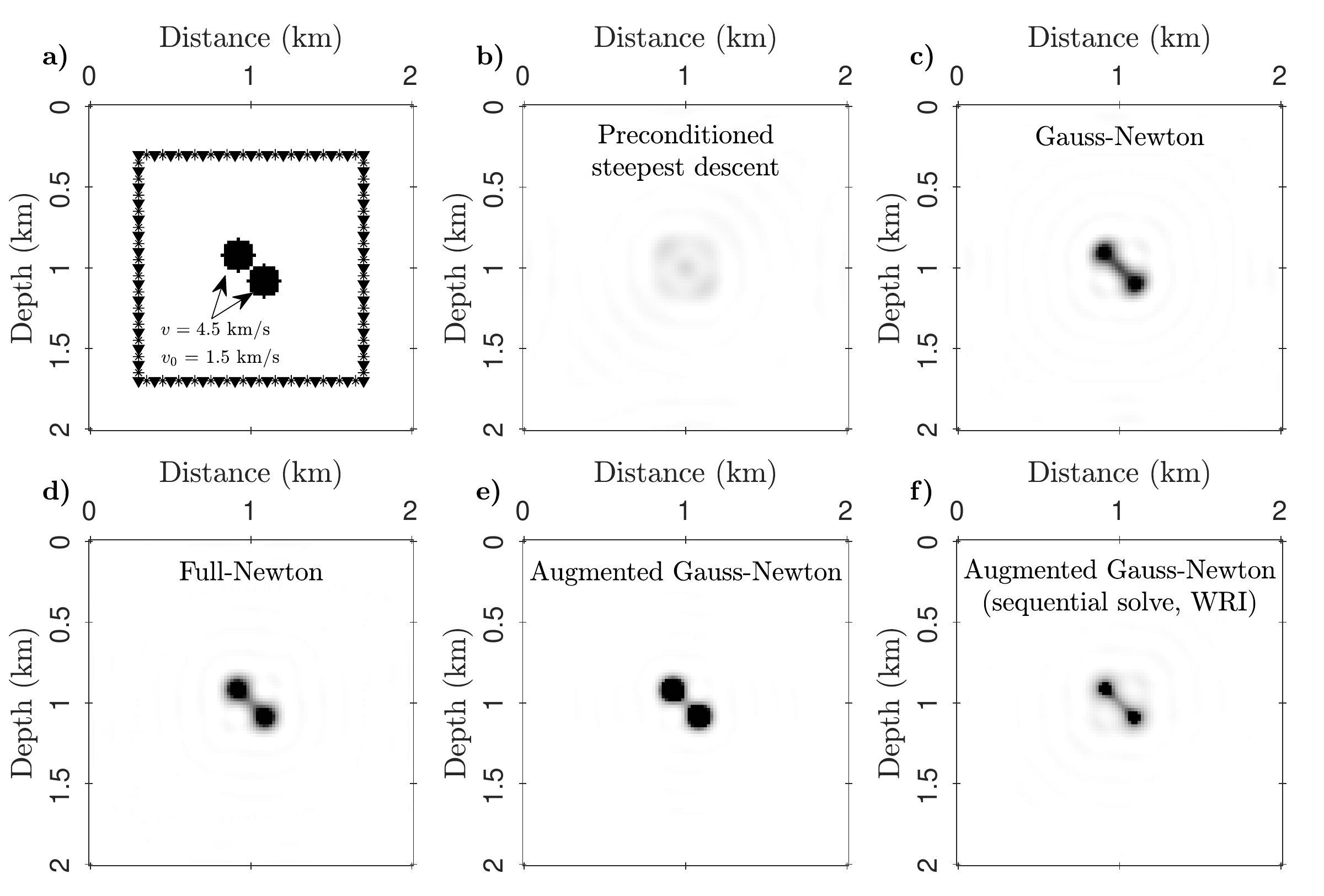}
\caption{(a) Inclusion model with sources (stars) and receivers (triangles) location. Estimated velocity models by (b) preconditioned steepest-descent, (c) Gauss-Newton, (d) full-Newton, (e) augmented Gauss-Newton, (f) augmented Gauss-Newton method with sequential solve or WRI.}
\label{fig:FIG2}
\end{figure}
\begin{figure}[htb!]
\center
\includegraphics[width=0.5\columnwidth]{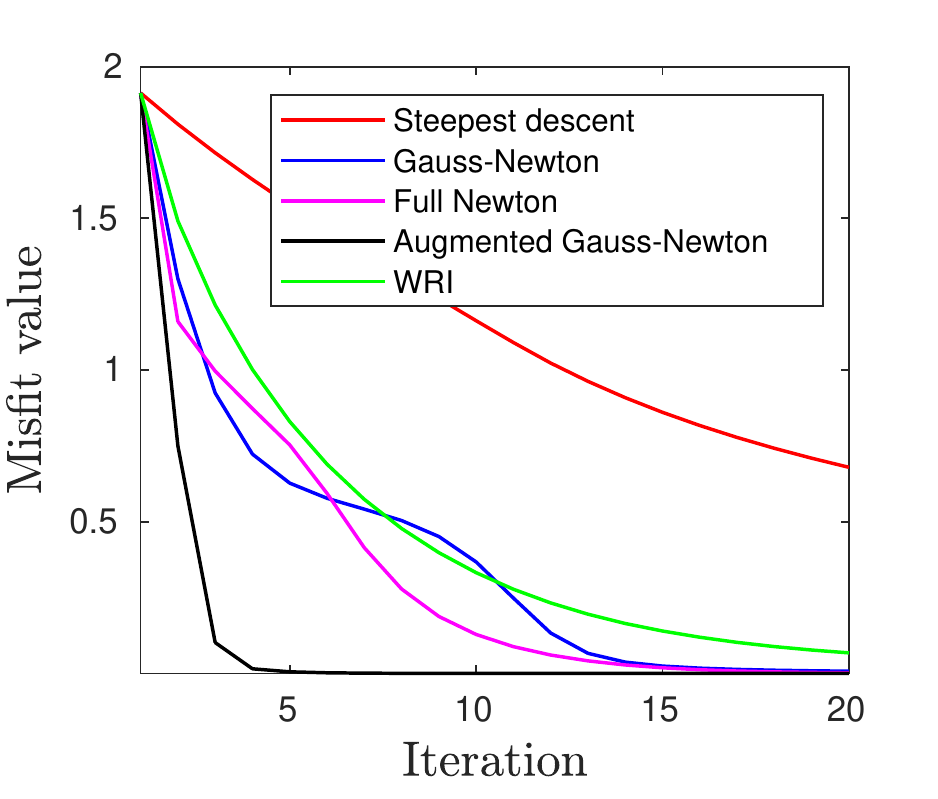}
\caption{Evolution of the misfit function for different methods. Preconditioned steepest-descent method (red line), Gauss-Newton method (blue line), Full-Newton method (magenta line), augmented Gauss-Newton method (black line), and augmented Gauss-Newton method with sequential solve or WRI (green line).}
\label{fig:FIG3}
\end{figure}
%
\section{Numerical results}
We consider the estimation of the velocity model shown in \fref{fig:FIG2}a using frequency-space FWI algorithm with different approximate Hessian matrices. 
The subsurface model contains two circular anomalies of velocity 4.5 km/s embedded in a homogeneous background of velocity 1.5 km/s (\fref{fig:FIG2}a). 
The dimensions of the model are 2.0 km in distance and 2.0 km in depth, and the grid spacing is 20 m.
The acquisition consists of 112 equally-spaced sources positioned around the model (shown by stars in \fref{fig:FIG2}a); and the wavefield due to each source is recorded by 112 equally-spaced receivers around the model (shown by triangles in \fref{fig:FIG2}a). 

We start the inversion from the homogeneous background model and perform 20 FWI iterations over one dataset corresponding to a 5 Hz frequency.
The simplest algorithm, based upon preconditioned steepest descent, uses the pseudo-Hessian approximation $\bold{H}\approx\bold{L}^T\bold{L}$ to update the model with \eref{NewtonSystem}. The associated result shows that this algorithm failed to correctly identify the two inclusion anomalies (\fref{fig:FIG2}b). 
We also use the GN method, which is based on the more advanced Hessian matrix $\bold{H}\approx\bold{L}^T\bold{S}^T\bold{SL}$, introducing the first-order derivative information of the misfit function. The associated result (\fref{fig:FIG2}c) shows significant improvement gained in the quality of the model estimate compared with that of the steepest-descent method, \fref{fig:FIG2}b. 
Incorporating the second-order derivative information is helpful in correct estimation of the large contrast anomalies presented in this model \citep{Metivier_2012_SOAS}.
Thus, we apply the full-Newton method by using the exact Hessian matrix and the resulting model is shown in \fref{fig:FIG2}d. With the help of the nonlinear term of the Hessian matrix, the two anomalies have been well resolved. The improvement achieved by using more accurate Hessian matrices are also highlighted by the convergence curves shown is \fref{fig:FIG3}.

We continued by examining the performance of the proposed augmented Hessian.
In order to show the effect of sequential solve, we use the AGN matrix in \eref{EGN_app}, which uses only the $\bold{R}^{11}$ part of the nonlinear term: $\bold{H} \approx   \bold{L}^T \bold{S}^T\bold{S}\bold{L} + \bold{L}^T \bold{S}^T\bold{S}\delta\bold{L}$. The result obtained with this augmented Hessian shows significant improvement (\fref{fig:FIG2}e) and is even better than that of the full-Newton method, \fref{fig:FIG2}d, which may be explained by numerical instability of the full Hessian matrix.  
Applying the sequential solve procedure gives the same updating formula as WRI, \eref{GN_update_seq_data} or \ref{WRIm3}. The result obtained by the sequential solve is shown in \fref{fig:FIG2}f. By comparing Figures \ref{fig:FIG2}e and f, we see that the quality of the estimate is decreased due to the approximate solve of the Newton's system obtained by the sequential procedure. This is also evident from the associated convergence curves in \fref{fig:FIG3}. 
However, the sequential solve procedure is computationally more efficient compared with the direct solve. More in depth comparison between these approaches remain to be done by using more challenging FWI problems in the presence of inaccurate initial models, where the nonlinear term of the Hessian play significant role for convergence to accurate solutions.

We also tested the performance of the FWI with the AGN Hessian versus WRI by using a near-surface imaging example \citep{Metivier_2017_TRU}. The subsurface model is composed of a homogeneous background of velocity 300 m/s and
two superimposed concrete structures of velocity 4000 m/s (\fref{fig:FIG3}a). 
The dimensions of the model are 15 m in distance and 3 m in depth, and the grid spacing is 0.15 m.
The acquisition consists of 120 equally-spaced sources/receivers positioned around the model. 
We start the inversion from the homogeneous background model and perform 200 FWI iterations over data sets corresponding to the frequencies 100, 125, 150, 175, 200, 225, 250, 275, and 300 Hz. 
The result obtained by FWI with the AGN Hessian and WRI are shown in Fig. \ref{fig:FIG3}b and c, respectively. The associated convergence curves are compared in Fig. \ref{fig:FIG3}d showing better performance of the  FWI with the AGN Hessian. 
\begin{figure}[htb!]
\center
\includegraphics[width=1\columnwidth]{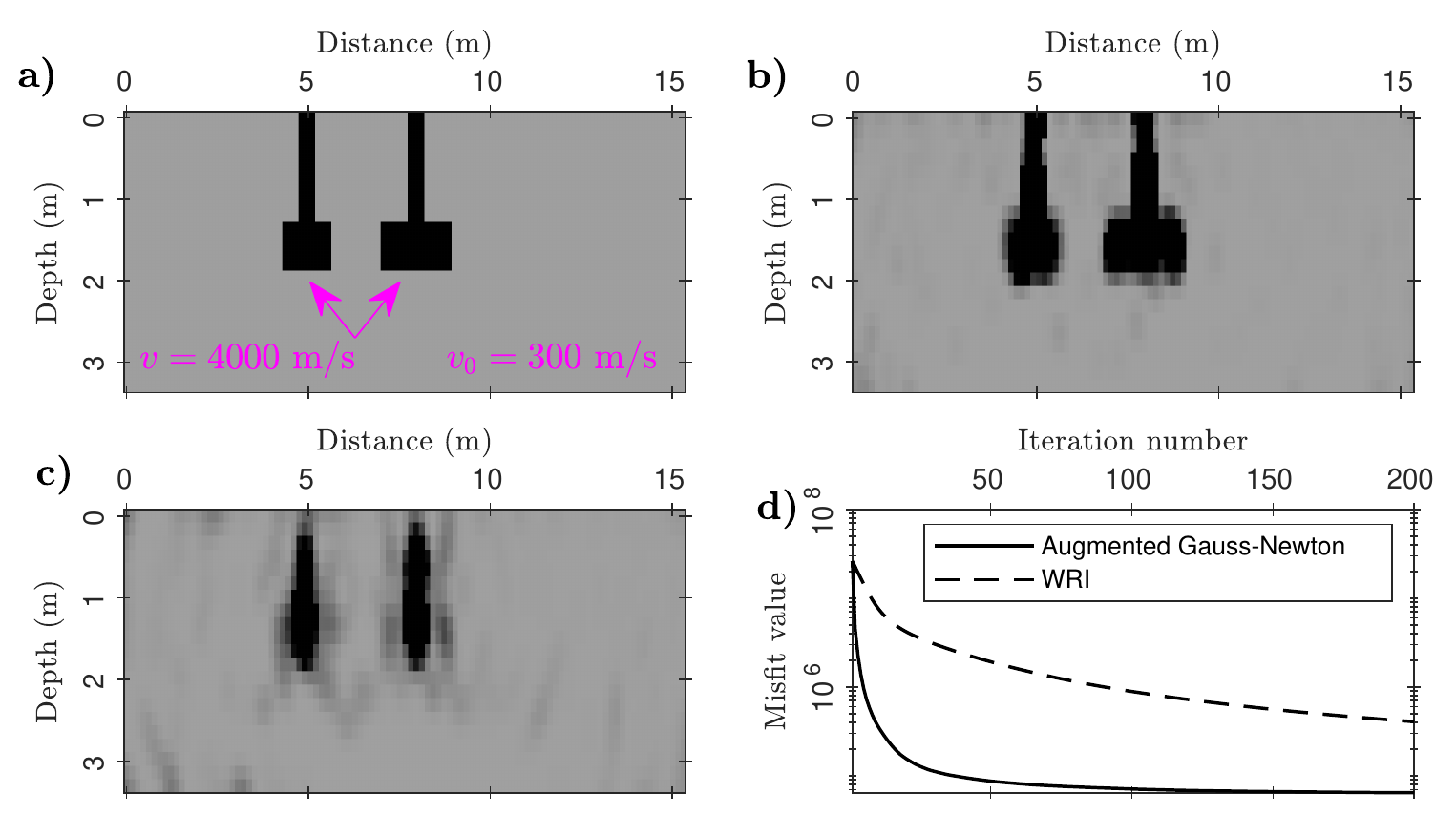}
\caption{(a) True near-surface model. Estimated velocity models by (b) augmented Gauss-Newton method and (c) augmented Gauss-Newton method with sequential solve or WRI. (d) Evolution of the misfit function versus iteration.}
\label{fig:FIG3}
\end{figure}
%
%
\section{Conclusion}
We proposed a new approximation of the Hessian matrix for FWI as an augmentation of the standard Gauss-Newton (GN) Hessian including the second-order derivative information.
We have shown numerically that this new Hessian is more well-posed compared with the original Hessian, and thus leads to more stable second-order updating formula.
Moreover, a special form of the proposed Newton system, when solved approximately, leads to an update formula which is exactly the same as the update formula provided by the WRI. 
Future research will investigate the characteristics of different terms of the New Hessian and developing efficient algorithms for solving the associated Newton system.
%

\textbf{Acknowledgments:} This study was funded by the WIND consortium (\hyperlink{https://www.geoazur.fr/WIND}{www.geoazur.fr/WIND}). 

\newcommand{\SortNoop}[1]{}

\end{document}